# EULER HYDRODYNAMICS OF ONE-DIMENSIONAL ATTRACTIVE PARTICLE SYSTEMS

By C. Bahadoran, H. Guiol, K. Ravishankar and E. Saada

*Université Clermont-Ferrand 2, INP Grenoble, SUNY New Paltz and CNRS Rouen*

We consider attractive irreducible conservative particle systems on $\mathbb{Z}$, without necessarily nearest-neighbor jumps or explicit invariant measures. We prove that for such systems, the hydrodynamic limit under Euler time scaling exists and is given by the entropy solution to some scalar conservation law with Lipschitz-continuous flux. Our approach is a generalization of Bahadoran et al. [*Stochastic Process. Appl.* **99** (2002) 1–30], from which we relax the assumption that the process has explicit invariant measures.

**1. Introduction.** In this paper we study the hydrodynamic behavior of a class of asymmetric particle systems of $\mathbb{Z}$, which arise as a natural generalization of the asymmetric exclusion process. For the latter a variety of results and approaches are available: the hydrodynamic limit is given by the entropy solutions to the scalar conservation law

$$(1) \qquad \partial_t u(x,t) + \partial_x G(u(x,t)) = 0,$$

where $u(\cdot,\cdot)$ is the density field, $G$—the *macroscopic flux*—is given by $G(u) = \gamma u(1-u)$, and $\gamma$ is the mean drift of a particle. Because there is a single conserved quantity (i.e., mass) for the particle system, and an ergodic equilibrium measure for each density value, (1) can be guessed through heuristic arguments if one takes for granted that the system is in *local equilibrium*. The macroscopic flux $G$ is obtained by an equilibrium expectation of a microscopic flux which can be written down explicitly from the dynamics. For the simple exclusion process, equilibrium measures are also explicit, and so is $G$. A rigorous proof of the hydrodynamic limit turns out to be a difficult problem, mainly because of the nonexistence of strong solutions for (1) and the nonuniqueness of weak solutions. Since the conservation law is not sufficient to pick a single solution, the so-called entropy weak solution must be characterized by additional properties; one must then look for related properties of the particle system to establish its convergence to the entropy solution.

The first approach to this problem was developed among others by Andjel and Vares [2], who proved the hydrodynamic limit (1) for the asymmetric exclusion process when the initial datum is a single step with density $\lambda$ to the left and $\rho$ to the right. This is usually referred to as the *Riemann problem* for (1), in which case the solution has a simple explicit form when $G$ is strictly convex or concave. Their approach combined this simplifying feature and comparison with equilibrium systems, which are explicitly known. Comparison is allowed by the monotonicity (also called *attractiveness*) of the system. The method of Andjel and Vares [2] also applies to other attractive systems with explicit invariant measures, such as the zero-range process; however in this case, one must add the assumption of flux concavity or convexity, which is not systematically true, nor easy to verify. Note that the $k$-step exclusion process introduced by Guiol [12] is an attractive process with product invariant measures but simple nonconcave and nonconvex flux function.

A breakthrough was achieved by Rezakhanlou [20], who proved the hydrodynamic limit for the asymmetric exclusion process in any dimension, with an arbitrary initial datum in (1). His approach is more general and based on entropy inequalities. It also applies, without any convexity assumptions, to other attractive processes with product invariant measures.

In [5], we initiated a "resurrection" of the approach of Andjel and Vares [2]. Our goal was to show that, with a little more work, the same result could be reached in dimension one with the same approach as in [20], but in a constructive way. Our argument was based on two ideas. First, we relaxed the convexity assumption of Andjel and Vares [2] by introducing a variational formula for the Riemann solution. Next, we provided a general argument, for finite-range attractive processes, showing that the hydrodynamic limit for Riemann initial condition implied hydrodynamics for general initial datum. This argument is inspired by the Glimm's scheme in the theory of hyperbolic conservation laws (see, e.g., Chapter 5 of [25]), a procedure to reconstruct general entropy solutions from Riemann entropy solutions.

In the present paper we refine our approach to obtain new results. Namely we establish the hydrodynamic limit for quite general finite-range attractive processes whose invariant measures are not product, nor even explicit at all, and not much is known about their properties: for instance the support of the densities of extremal elements is unknown. This extends two former results in this direction: the first one by Seppäläinen [23] for the totally asymmetric, nearest-neighbor $K$-exclusion process, based on the author's variational coupling method; the second one by Rezakhanlou [21] for the partially asymmetric nearest-neighbor $K$-exclusion process, using an abstract characterization of Hamilton–Jacobi semigroups, and the ergodic theorem of Ekhaus and Gray [9] for nearest-neighbor attractive systems.

In a forthcoming paper, we will demonstrate how a further refinement of our method leads to a strong law of large numbers (i.e., almost sure hydrodynamics) for attractive systems on $\mathbb{Z}$.

**2. Notation and results.** Throughout this paper $\mathbb{N} = \{1, 2, \ldots\}$ will denote the set of natural numbers, $\mathbb{Z}^+ = \{0, 1, 2, \ldots\}$ the set of nonnegative integers, and $\mathbb{R}^{+*} = \mathbb{R}^+ - \{0\}$ the set of positive real numbers. We consider particle systems on $\mathbb{Z}$ with at most $K$ particles per site, $K \in \mathbb{N}$. Thus the state space of the process is $\mathbf{X} = \{0, 1, \ldots, K\}^{\mathbb{Z}}$, which we endow with the product topology. A function defined on $\mathbf{X}$ is called *local* if it depends on the variable $\eta \in \mathbf{X}$ only through $(\eta(x), x \in \Lambda)$ for some finite subset $\Lambda$ of $\mathbb{Z}^d$. If $\eta$ is an $\mathbf{X}$-valued random variable and $\nu$ a probability measure on $\mathbf{X}$, we write $\eta \sim \nu$ to specify that $\eta$ has distribution $\nu$. The notation $\nu(f)$, where $f$ is a real-valued function and $\nu$ a probability measure on $\mathbf{X}$, will be an alternative for $\int_{\mathbf{X}} f \, d\nu$. We write $\nu_n \Rightarrow \nu$ to denote weak convergence of a sequence $(\nu_n, n \in \mathbb{N})$ of probability measures on $\mathbf{X}$ to some probability measure $\nu$ on $\mathbf{X}$, that is, $\nu_n(f) \to \nu(f)$ as $n \to \infty$ for every continuous function $f$ on $\mathbf{X}$.

2.1. *The model.* The dynamics consists of particles' jumps, according to the Markov generator

$$(2) \qquad Lf(\eta) = \sum_{x,y \in \mathbb{Z}} p(y-x) b(\eta(x), \eta(y)) [f(\eta^{x,y}) - f(\eta)]$$

for a local function $f$, where $\eta^{x,y}$ denotes the new state after a particle has jumped from $x$ to $y$ [i.e., $\eta^{x,y}(x) = \eta(x) - 1$, $\eta^{x,y}(y) = \eta(y) + 1$, $\eta^{x,y}(z) = \eta(z)$ otherwise], $p$ is the particles' jump kernel, that is, $\sum_{z \in \mathbb{Z}} p(z) = 1$, and $b : \mathbb{Z}^+ \times \mathbb{Z}^+ \to \mathbb{R}^+$ is the jump rate. We assume that $p$ and $b$ satisfy:

(A1) The semigroup of $\mathbb{Z}$ generated by the support of $p$ is $\mathbb{Z}$ itself (standard irreducibility);

(A2) $p$ is finite range, that is, there exists $M > 0$ such that $p(x) = 0$ for all $|x| > M$;

(A3) $b(0, \cdot) = 0$, $b(\cdot, K) = 0$ (no more than $K$ particles per site), and $b(i,j) > 0$ for $0 < i \leq K$ and $0 \leq j < K$ (nondegenerate jump rates);

(A4) $b$ is nondecreasing (nonincreasing) in its first (second) argument.

REMARK. In view of assumption (A4), the third condition in (A3) can be equivalently replaced by the simpler condition $b(1, K-1) > 0$.

If $b$ satisfies additional algebraic relations (see [8]), which include, for example, the asymmetric exclusion process, the system has an explicit family of invariant product measures indexed by the mean density of particles, and is the so-called *misanthrope's process*. Here we do not assume these relations, and the system in general has no explicit invariant measures. This is even the case for the simplest possible $b$, that is, $b(n,m) = \mathbf{1}_{\{n>0, m<K\}}$, for which the system is called $K$-*exclusion process*—see [23]—(except for symmetric $p$, in which case there are still product invariant measures, see [13]; but then the relevant time scale is diffusive and not Euler, and Theorem 2.2 below yields a trivial nonevolving hydrodynamic limit).

The special form

$$p(y-x)b(\eta(x), \eta(y)) \tag{3}$$

in (2), in which a jump kernel $p(\cdot)$ and a configuration-dependent part $b(\cdot,\cdot)$ are decoupled, does not play any role in our paper, but we retained it by analogy with [8], where it was needed to exhibit product invariant measures for a certain class of $b$'s. Here we might as well consider more general jump rates of the form

$$b(y-x, \eta(x), \eta(y)) \tag{4}$$

for which the set of assumptions (A1)–(A4) would have to be replaced by:

(A1$'$) The semigroup of $\mathbb{Z}$ generated by the set $\{z \in \mathbb{Z} : \inf_{n>0, m<K} b(z, n, m) > 0\}$ is $\mathbb{Z}$ itself.

(A2$'$) Finite range assumption: there exists $M > 0$ such that $b(z, \cdot, \cdot) = 0$ for all $|z| > M$.

(A3$'$) $b(\cdot, 0, \cdot) = 0$, $b(\cdot, \cdot, K) = 0$ (no more than $K$ particles per site).

(A4$'$) $b$ is nondecreasing (nonincreasing) in its second (third) argument.

Note that, in this more general context, assumption (A1$'$) above replaces both assumption (A1) and the second part of assumption (A3). Under these assumptions, the coupling arguments of Cocozza-Thivent [8] needed in Section 3 below carry over. There is another reason why (3) is not so relevant here: for the misanthrope's process of Cocozza-Thivent [8], the decoupling in (3) is reflected in the hydrodynamic limit in the form of a decoupling in the macroscopic flux, that is

$$G(u) = \gamma H(u), \tag{5}$$

where $\gamma := \sum_z z p(z)$ is the mean drift depending only on $p(\cdot)$, and $H(\cdot)$ depends only on $b(\cdot)$. This is in fact due to existence of product invariant measures for all density values, see Remark 1. As a result, $p(\cdot)$ and $b(\cdot,\cdot)$ have separate interpretations in the hydrodynamic limit. In our case, without product invariant measures, there is no a priori reason for a decoupling like

(5) in the macroscopic flux, hence no particular sense in the microscopic decoupling (3).

Let the coordinatewise partial order on **X** be defined by $\eta \leq \xi$ if and only if $\eta(x) \leq \xi(x)$ for every $x \in \mathbb{Z}$. This induces a partial stochastic order for probability measures $\mu_1$ and $\mu_2$ on **X**; namely, we write $\mu_1 \leq \mu_2$ if the following equivalent conditions hold (see, e.g., [19]):

(i) For every nondecreasing nonnegative function $f$ on **X**, $\mu_1(f) \leq \mu_2(f)$.

(ii) There exists a coupling measure $\widetilde{\mu}$ on the product space $\widetilde{\mathbf{X}} = \mathbf{X} \times \mathbf{X}$ with marginals $\mu_1$ and $\mu_2$, such that $\widetilde{\mu}\{(\eta,\xi): \eta \leq \xi\} = 1$.

It is shown in [8] that assumption (A4) implies the existence of a coupled Markov generator $\tilde{L}$ (called the *basic coupling*) on $\widetilde{\mathbf{X}}$ such that: (a) $\tilde{L}$ generates Feller processes $(\eta_t, \xi_t)_{t \geq 0}$ whose components are Feller processes generated by $L$, and (b) for these coupled processes, $\eta_0 \leq \xi_0$ a.s. implies $\eta_t \leq \xi_t$ a.s. for every $t > 0$. Properties (a) and (b) above imply *monotonicity* of the semigroup $S(\cdot)$, that is,

(6) $$\mu \leq \nu \implies \forall t > 0, \mu S(t) \leq \nu S(t).$$

Either properties (a) and (b), or (6), are usually called *attractiveness*. These are still true for the more general model (4) under assumption (A4').

2.2. *Scalar conservation laws and entropy solutions.* We recall the definition of entropy solutions to scalar conservation laws, which will appear as hydrodynamic limits of the above models. For more details, we refer to the textbooks of Godlewski and Raviart [11], Serre [25] or Bressan [6].

Let $G: [0,K] \to \mathbb{R}$ be a Lipschitz-continuous function, called the *flux*. It is a.e. differentiable, and its derivative $G'$ is an (essentially) uniformly bounded function. We consider the scalar conservation law

(7) $$\partial_t u + \partial_x [G(u)] = 0,$$

where $u = u(x,t)$ is some $[0,K]$-valued density field defined on $\mathbb{R} \times \mathbb{R}^+$. Equation (7) has no strong solutions in general: even starting from a smooth Cauchy datum $u(\cdot, 0) = u_0$, discontinuities (called shocks in this context) appear in finite time. Therefore it is necessary to consider weak solutions, but then uniqueness is lost for the Cauchy problem. To recover uniqueness, we need to define *entropy solutions*.

Let $\phi: [0,K] \to \mathbb{R}$ be a convex function. In the context of hyperbolic systems, such a function is called an *entropy*. We define the associated *entropy flux* $\psi$ on $[0,K]$ as

$$\psi(u) := \int_0^u \phi'(v) G'(v) \, dv,$$

$(\phi, \psi)$ is called an *entropy-flux pair*. A Borel function $u : \mathbb{R} \times \mathbb{R}^{+*} \to [0, K]$ is called an *entropy solution* to (7) if and only if it is entropy-dissipative, that is,

$$\partial_t \phi(u) + \partial_x \psi(u) \leq 0 \tag{8}$$

in the sense of distributions on $\mathbb{R} \times \mathbb{R}^{+*}$ for any entropy-flux pair $(\phi, \psi)$. Note that, by taking $\phi(u) = \pm u$ and hence $\psi(u) = \pm G(u)$, we see that an entropy solution is indeed a weak solution to (7). This definition can be motivated by the following points: (i) when $G$ and $\phi$ are continuously differentiable, (7) implies equality in strong sense in (8) (this follows from the chain rule for differentiation); (ii) this no longer holds in general if $u$ is only a weak solution to (7); (iii) the inequality (8) can be seen as a macroscopic version of the second law of thermodynamics that selects physically relevant solutions. Indeed, one should think of the *concave* function $h = -\phi$ as a thermodynamic entropy, and spatial integration of (8) shows that the total thermodynamic entropy may not decrease during the evolution (this is rigorously true for periodic boundary conditions, in which case the total entropy is well defined).

Kružkov proved the following fundamental existence (Theorem 2 of [15]) and uniqueness (Theorem 5 of [15]) result:

THEOREM 2.1. *Let $u_0 : \mathbb{R} \to [0, K]$ be a Borel measurable initial datum. Then there exists a unique (up to a Lebesgue-null subset of $\mathbb{R} \times \mathbb{R}^{+*}$) entropy solution $u$ to (7) subject to the initial condition*

$$\lim_{t \to 0^+} u(\cdot, t) = u_0(\cdot) \qquad \text{in } L^1_{\text{loc}}(\mathbb{R}). \tag{9}$$

*This solution [has a representative in its $L^\infty(\mathbb{R} \times \mathbb{R}^{+*})$ equivalence class that] is continuous as a mapping $t \mapsto u(\cdot, t)$ from $\mathbb{R}^{+*}$ to $L^1_{\text{loc}}(\mathbb{R})$.*

We recall here that a sequence $(u_n, n \in \mathbb{N})$ of Borel measurable functions on $\mathbb{R}$ is said to converge to $u$ in $L^1_{\text{loc}}(\mathbb{R})$ if and only if

$$\lim_{n \to \infty} \int_I |u_n(x) - u(x)| \, dx = 0$$

for every bounded interval $I \subset \mathbb{R}$.

REMARK. Kružkov's theorems are stated for a continuously differentiable $G$. However the proof of the uniqueness result (Theorem 2 of [15]) uses only Lipschitz continuity. In the Lipschitz-continuous case, existence could be derived from Kružkov's result by a flux approximation argument. However a different, self-contained (and constructive) proof of existence in this case can be found in Chapter 6 of [6].

The uniqueness statement of Theorem 2.1 is in fact a consequence of the following $L^1$-*contraction* and *finite propagation* properties (see also property (2) in Theorem 3.1 of [5], or [11] or [25]):

PROPOSITION 2.1. *Let $V := \|G'\|_\infty$ denote the Lipschitz constant of the macroscopic flux $G(\cdot)$ in (7). Then, for any $x < y$ and $0 < t < (y-x)/(2V)$,*

(10) $$\int_{x+Vt}^{y-Vt} |u(z,t) - v(z,t)|\,dz \leq \int_x^y |u_0(z) - v_0(z)|\,dz,$$

*where $u(\cdot,\cdot)$ and $v(\cdot,\cdot)$, respectively, denote the entropy solutions to (7) with initial data $u_0(\cdot)$ and $v_0(\cdot)$. In particular, if $u_0(\cdot)$ and $v_0(\cdot)$ coincide a.e. on $[x,y]$, then $u(\cdot,t)$ and $v(\cdot,t)$ coincide a.e. on $[x+Vt, y-Vt]$.*

We next recall a possibly more familiar definition of entropy solutions based on shock admissibility conditions, but valid only for solutions with bounded variation. This point of view selects the relevant weak solutions by specifying what kind of discontinuities are permitted. The following geometric condition is known as Oleĭnik's entropy condition (see, e.g., [11] or [25]). A discontinuity $(u^-, u^+)$, with $u^\pm := u(x \pm 0, t)$, is called an entropy shock, if and only if:

(11) The chord of the graph of $G$ between $u^-$ and $u^+$ lies: below the graph if $u^- < u^+$, above the graph if $u^+ < u^-$.

In the above condition, "below" or "above" is meant in wide sense, that is does not exclude that the graph and chord coincide at some points between $u^-$ and $u^+$. In particular, when $G$ is strictly convex (resp. concave), one recovers the fact that only (and all) decreasing (resp. increasing) jumps are admitted. Note that, if the graph of $G$ is linear on some nontrivial interval, condition (11) implies that any increasing or decreasing jump within this interval is an entropy shock.

Condition (11) can be used to select entropy solutions among weak solutions. Let $\mathrm{TV}_I$ denote the variation of a function defined on some closed interval $I = [a,b] \subset \mathbb{R}$, that is,

$$\mathrm{TV}_I[u(\cdot)] = \sup_{x_0=a<x_1<\cdots<x_n=b} \sum_{i=0}^{n-1} |u(x_{i+1}) - u(x_i)|.$$

Let us say that $u = u(\cdot,\cdot)$ defined on $\mathbb{R} \times \mathbb{R}^{+*}$ has locally bounded space variation if

(12) $$\sup_{t \in J} \mathrm{TV}_I[u(\cdot, t)] < +\infty$$

for every bounded closed space interval $I \subset \mathbb{R}$ and bounded time interval $J \subset \mathbb{R}^{+*}$. Then the following result, which will be used in Section 4.1, is a consequence of Vol'pert [26].

PROPOSITION 2.2. *Let $u$ be a weak solution to (7) with locally bounded space variation. Then $u$ is an entropy solution to (7) if and only if, for a.e. $t > 0$, all discontinuities of $u(\cdot, t)$ are entropy shocks.*

For completeness a proof of this statement is given in Appendix A.

REMARK. One can show that, if the Cauchy datum $u_0$ has locally bounded variation, the unique entropy solution given by Theorem 2.1 has locally bounded space variation. Hence Proposition 2.2 extends into an existence and uniqueness theorem within functions of locally bounded space variation, where entropy solutions may be defined as weak solutions satisfying (11), without reference to (8).

2.3. *The hydrodynamic limit.* Before stating our main result, we recall some standard definitions in hydrodynamic limits. The integer $N \in \mathbb{N}$ is the scaling parameter for the hydrodynamic limit, that is, the inverse of the macroscopic distance between two consecutive sites. The empirical measure of a configuration $\eta$ viewed on scale $N$ is given by

$$(13) \qquad \alpha^N(\eta, dx) = N^{-1} \sum_{y \in \mathbb{Z}} \eta(y) \delta_{y/N}(dx) \in \mathcal{M},$$

where $\mathcal{M}$ denotes the set of positive, locally finite measures on $\mathbb{R}$. The set $\mathcal{M}$ is equipped with the topology of vague convergence, defined by convergence for continuous test functions with compact support. Let $u(\cdot)$ be a deterministic bounded Borel function on $\mathbb{R}$. A sequence of random configurations $\eta^N$ is said to have *density profile* $u(\cdot)$ if $\alpha^N(\eta^N, dx)$ converges in probability to $u(\cdot) dx$ as $N \to \infty$.

Let $u(\cdot, \cdot)$ be a deterministic bounded Borel function on $\mathbb{R} \times \mathbb{R}^+$. A sequence of processes $\eta^N_\cdot = (\eta^N_t, t \geq 0)$ generated by $L$, with random initial configurations $\eta^N_0$, is said to have *hydrodynamic limit* $u(\cdot, \cdot)$ (under Euler time scaling), if for every $t \geq 0$, the sequence $\eta^N_{Nt}$ has density profile $u(\cdot, t)$.

We now state our main result.

THEOREM 2.2. *Assume the sequence $\eta^N_0$ has density profile $u_0(\cdot)$, where $u_0(\cdot)$ is a measurable $[0, K]$-valued profile on $\mathbb{R}$. Then the sequence of processes $\eta^N_\cdot$ has a hydrodynamic limit $u(\cdot, \cdot)$ given by the unique entropy solution to the scalar conservation law (7) with initial condition $u_0$, where $G : [0, K] \to \mathbb{R}^+$ is a Lipschitz-continuous flux function defined below.*

The flux function $G$ in Theorem 2.2 is obtained as follows. Let $\mathcal{I}$ and $\mathcal{S}$ denote respectively the set of invariant measures for $L$, and the set of shift-invariant measures on $\mathbf{X}$. Since $\mathbf{X}$ is compact, and the process has the Feller property, $\mathcal{I} \cap \mathcal{S}$ is a nonempty (the empty and full configurations are unchanged by the dynamics) convex compact subset of the set of probability measures on **X**. Thus any element $\nu \in \mathcal{I} \cap \mathcal{S}$ is a mixture of its extremal elements. In the usual setting of Andjel and Vares [2], Rezakhanlou [20], Bahadoran et al. [5], it is known that

$$(\mathcal{I} \cap \mathcal{S})_e = \{\nu_\rho, \rho \in [0, K]\},$$

where $\nu_\rho$ are product measures and the index $e$ denotes the extremal elements. In our setting, there are no explicit invariant measures. However, as we will explain in Section 3, the same coupling arguments as used by Liggett [18], relying only on attractiveness and irreducibility, are still sufficient to establish

$$(14) \qquad (\mathcal{I} \cap \mathcal{S})_e = \{\nu_\rho, \rho \in \mathcal{R}\},$$

where $\mathcal{R}$ is a closed subset of $[0, K]$ containing 0 and $K$, and $\nu_\rho$ is a shift-invariant probability measure on **X** with $\nu_\rho[\eta(0)] = \rho$. Whether $\mathcal{R} = [0, K]$ is an open problem.

We define the *microscopic* flux (here across site 0) by

$$(15) \quad j(\eta) = \sum_{x \leq 0, y > 0} p(y-x)b(\eta(x), \eta(y)) - \sum_{x > 0, y \leq 0} p(y-x)b(\eta(x), \eta(y)).$$

The term "microscopic flux" is justified by the microscopic conservation equation [clear from (2)]

$$(16) \qquad L\left[\sum_{z=x+1}^{y} \eta(z)\right] = \tau_x j(\eta) - \tau_y j(\eta)$$

for $x < y$, where $\tau$ is the space shift, that is, $\tau_x \eta(z) = \eta(x+z)$ and $\tau_x j(\eta) = j(\tau_x \eta)$. Equation (16) says that the instantaneous variation of the number of particles in interval $(x, y] \cap \mathbb{Z}$ is due to the fluxes (translates of $j$) at the boundaries. The *macroscopic flux* is now defined by: for $\rho \in \mathcal{R}$,

$$(17) \qquad G(\rho) = \nu_\rho[j(\eta)],$$

then interpolate $G$ linearly on the complement of $\mathcal{R}$, which is an at most countable union of disjoint open intervals. Note that, for $\rho \in \mathcal{R}$, by shift-invariance of $\nu_\rho$, we also have $G(\rho) = \nu_\rho[j'(\eta)]$ with $j'(\eta) = \sum_z z p(z) b(\eta(0), \eta(z))$.

REMARK 1. For the misanthrope's process, the decoupled form (5) follows from (15), (17), and the fact that $\nu_\rho$ exists and is a product measure for all $\rho$. Indeed, we can then define $H(\rho) := \nu_\rho[b(\eta(x), \eta(y))]$ for any $x \neq y$, because this quantity is independent of $x$ and $y$. Hence (5) holds with

$$\gamma := \sum_{x \leq 0, y > 0} p(y-x) - \sum_{x > 0, y \leq 0} p(y-x) = \sum_z z p(z).$$

REMARK 2. *Unbounded systems.* The arguments and results of Section 3 do not use the assumption that the number of particles per site is bounded by $K < +\infty$. Thus they are still true for systems with state space $\mathbb{N}^{\mathbb{Z}}$ satisfying assumptions (A1), (A2), (A4) and

(A3″) $b(0, \cdot) = 0$ and $b(i, j) > 0$ for $i > 0$.

In this case (14) still holds with $\mathcal{R}$ a closed subset of $\mathbb{R}^+$. However, in order to define a macroscopic flux by interpolation, we would need to know that $\mathcal{R}$ is unbounded, which we cannot establish. Then the arguments and results of Section 4 would also extend, thus establishing Riemann hydrodynamics. On the other hand, we do not currently know how to carry out the passage from Riemann to general hydrodynamics (Section 5) for unbounded systems without product invariant measures, because we cannot prove Lemma 5.1 for such systems.

REMARK 3. Nothing can be said in general about the set $\mathcal{R}$, except that it is closed and contains 0 and $K$. However for the totally asymmetric nearest-neighbor $K$-exclusion process, one can benefit from some properties of $G$ derived in [23], to get slightly more precise information about $\mathcal{R}$.

COROLLARY 2.1. *For the totally asymmetric $K$-exclusion process, 0 and $K$ are limit points of $\mathcal{R}$, and $\mathcal{R}$ contains at least one point in $[1/3, K - 1/3]$.*

PROOF. The following properties were established in [23] for $G$:

(a) $G$ is symmetric around $u = K/2$: $G(K - u) = G(u)$ for every $u \in [0, K]$.

(b) $G$ is concave.

(c) $G$ has the following bounds: $F(u) \leq G(u) \leq H(u)$ for every $u \in [0, K]$, where

$$F(u) = \begin{cases} u(1-u), & \text{if } 0 \leq u \leq 1/2, \\ 1/4, & \text{if } 1/2 \leq u \leq K - 1/2, \\ (K-u)(1-(K-u)), & \text{if } K - 1/2 \leq u \leq K, \end{cases}$$

$$H(u) = \begin{cases} u/(1+u), & \text{if } 0 \leq u \leq K/2, \\ (K-u)/(1+K-u), & \text{if } K/2 \leq u \leq K. \end{cases}$$

If 0 were no limit point for $\mathcal{R}$, this would imply by construction of $G$ that there exists $\varepsilon > 0$ such that $G$ is linear on $[0, \varepsilon]$. But there can be no linear function lying between the lower bound $F$ and upper bound $H$ near the origin, because $F$ and $H$ have the same derivative at $u = 0$ and are strictly concave near the origin. The same argument works around $K$ by symmetry of $F$ and $H$.

Let us now assume $\mathcal{R} \cap [1/3, K - 1/3] = \varnothing$. Since $\mathcal{R}$ is closed, we must have $\mathcal{R} \cap [1/3 - \varepsilon, K - 1/3 + \varepsilon] = \varnothing$ for some $0 < \varepsilon < 1/3$. Then, by construction of $G$, the graph of $G$ between $u = 1/3 - \varepsilon$ and $u = K - 1/3 + \varepsilon$ must be a straight line. By the symmetry property (a) of $G$, this line must be horizontal, that is, $G$ has constant value over $[1/3 - \varepsilon, K - 1/3 + \varepsilon]$. Because of the lower bound in (c), this constant value must be at least $1/4$. Hence $G(1/3 - \varepsilon) \geq 1/4$, which contradicts the upper bound in (c), since $H(1/3 - \varepsilon) < H(1/3) = 1/4$. □

2.4. *Riemann solutions and Glimm's scheme.* Of special importance among entropy solutions are the solutions of the Riemann problem, that is, the Cauchy problem for particular initial data of the form

$$(18) \qquad u_0(x) = \lambda \mathbf{1}_{\{x<0\}} + \rho \mathbf{1}_{\{x \geq 0\}}.$$

Indeed: (i) these solutions can be computed explicitly and have a variational representation: see Section 4.1; (ii) one can construct approximations to the solution of the general Cauchy problem by using only Riemann solutions. This has inspired our belief that one could derive general hydrodynamics from Riemann hydrodynamics.

We will briefly explain here the principle of approximation schemes based on (ii), the most important of which is probably Glimm's scheme, introduced in [10]. Consider as initial datum a piecewise constant profile with finitely many jumps. The key observation is that, for small enough times, this can be viewed as a succession of noninteracting Riemann problems. To formalize this, we recall part of Lemma 3.4 of [5], whose proof was a direct application of Proposition 2.1 above. We denote by $R_{\lambda,\rho}(x,t)$ the entropy solution to the Riemann problem with initial datum (18).

LEMMA 2.1. *Let* $x_0 = -\infty < x_1 < \cdots < x_n < x_{n+1} = +\infty$, *and* $\varepsilon := \min_k(x_{k+1} - x_k)$. *Consider the Cauchy datum*

$$u_0 := \sum_{k=0}^n r_k \mathbf{1}_{(x_k, x_{k+1})},$$

*where* $r_k \in [0, K]$. *Then for* $t < \varepsilon/(2V)$, *with* $V$ *given in Proposition* 2.1, *the entropy solution* $u(\cdot, t)$ *at time* $t$ *coincides with* $R_{r_{k-1}, r_k}(\cdot - x_k, t)$ *on* $(x_{k-1} + Vt, x_{k+1} - Vt)$. *In particular,* $u(\cdot, t)$ *has constant value* $r_k$ *on* $(x_k + Vt, x_{k+1} - Vt)$.

Given some Cauchy datum $u_0$, we construct an approximate solution $\tilde{u}(\cdot, \cdot)$ for the corresponding entropy solution $u(\cdot, \cdot)$. To this end we define an approximation scheme based on a time discretization step $\Delta t > 0$ and a space discretization step $\Delta x > 0$. In the limit we let $\Delta x \to 0$ with the ratio $R := \Delta t / \Delta x$ kept constant under the condition

$$(19) \qquad\qquad\qquad R \leq 1/(2V)$$

known as the *Courant–Friedrichs–Lewy* (*CFL*) *condition*. Let $t_k := k\Delta t$ denote discretization times. We start with $k = 0$, setting $\tilde{u}_0^- := u_0$.

*Step* 1 (approximation step): Approximate $\tilde{u}_k^-$ with a piecewise constant profile $\tilde{u}_k^+$ whose step lengths are bounded below by $\Delta x$.

*Step* 2 (evolution step): For $t \in [t_k, t_{k+1})$, denote by $\tilde{u}_k(\cdot, t)$ the entropy solution at time $t$ with initial datum $\tilde{u}_k^+$ at time $t_k$. By (19) and Lemma 2.1, $\tilde{u}_k(\cdot, t)$ can be computed solving only Riemann problems. Set $\tilde{u}_{k+1}^- = \tilde{u}_k(\cdot, t_{k+1})$.

*Step* 3 (iteration): Increment $k$ and go back to step 1.

The approximate entropy solution is then defined by

$$\tilde{u}(\cdot, t) := \sum_{k \in \mathbb{N}} \tilde{u}_k(\cdot, t) \mathbf{1}_{[t_k, t_{k+1})}(t). \tag{20}$$

The efficiency of the scheme depends on how the approximation step is performed. In Glimm's scheme, the approximation $\tilde{u}_k^+$ is defined as

$$\tilde{u}_k^+ := \sum_{j \in k/2 + \mathbb{Z}} \tilde{u}_k^-((j + a_k/2)\Delta x) \mathbf{1}_{((j-1/2)\Delta x, (j+1/2)\Delta x)}, \tag{21}$$

where $a_k \in (-1, 1)$. Then we have the following convergence result.

THEOREM 2.3. *Let $u_0$ be a given measurable initial datum. Then every sequence $\varepsilon_n \downarrow 0$ as $n \to \infty$ has a subsequence $\delta_n \downarrow 0$ such that, for a.e. sequence $(a_k)$ w.r.t. product uniform measure on $(-1, 1)^{\mathbb{Z}^+}$, the Glimm approximation defined by (20) and (21) converges to $u$ in $L^1_{\mathrm{loc}}(\mathbb{R} \times \mathbb{R}^{+*})$ as $\Delta x = \delta_n \downarrow 0$.*

When $u_0$ has locally bounded variation, the above result is a specialization to scalar conservation laws of a more general result for systems of conservation laws: see Theorems 5.2.1, 5.2.2, 5.4.1 and comments following Theorem 5.2.2 in [25]. In Appendix B, we prove that it is enough to assume $u_0$ measurable.

Due to the nature of the approximation step (21), the proof of Theorem 2.3 does not proceed by direct estimation of the error between $\tilde{u}_k^\pm$ and $u(\cdot, t_k)$, but indirectly, by showing that limits of the scheme satisfy (8). Our approach is based on a different approximation, introduced first in Lemma 3.6 of [5], and refined in Lemma 5.5 of the present paper. This approximation allows direct control of the error by using a suitable distance, defined and denoted by $\Delta(u, v)$ in Section 5. Intuitively errors accumulate during approximation steps, but might be amplified by the evolution steps. The key properties of our approximation are that (i) the total error accumulated during the approximation step is negligible as $\varepsilon \to 0$; (ii) the error is not amplified by the evolution step, because the resolution semigroup of (7) is $\Delta$-contractant; (iii) a similar property holds at particle level, that is, the particle system is contractant for a microscopic version of $\Delta$. This allows us to mimic the scheme at particle level.

**3. Properties of invariant measures.** The main result of this section is

PROPOSITION 3.1. *(i) (14) holds with $\mathcal{R}$ a closed subset of $[0, K]$ containing $0$ and $K$, and $\nu_\rho$ a stationary shift-invariant measure such that $\nu_\rho[\eta(0)] = \rho$.*

*(ii) The measures $\nu_\rho$ are stochastically ordered: $\nu_\rho \leq \nu_{\rho'}$ if $\rho \leq \rho'$.*

*(iii) $\nu_\rho$ has a.s. density $\rho$, that is, $(2l+1)^{-1} \sum_{x=-l}^{l} \eta(x)$ converges $\nu_\rho$-a.s. to $\rho$ as $l \to \infty$.*

REMARK. A proof of (i) and (ii) can be found in [9] and [21] for nearest-neighbor processes. It is extended [as well as (iii)] to finite-range $K$-exclusion processes in [24]. In these papers it appears as a byproduct of a large time convergence result, which also implies an ergodic theorem. However if one is only interested in (i) and (ii), and not in the ergodic theorem, a much easier independent proof is possible, along the lines of [18] or [19], without using any convergence result. This is what we explain below.

If the elements of $(\mathcal{I} \cap \mathcal{S})_e$ were product measures, it would immediately follow that they are ergodic with respect to translations, implying (iii). Since this is not the case here, for a proof of (iii) we refer the reader to [21]; though in that paper the proof is given in the context of nearest-neighbor $K$-exclusion process, the same argument shows that, for any conservative finite-range system, every element of $(\mathcal{I} \cap \mathcal{S})_e$ has an a.s. density.

We divide the proof of (i)–(ii) into two lemmas.

LEMMA 3.1. *Any two elements $\mu_1$ and $\mu_2$ of $(\mathcal{I} \cap \mathcal{S})_e$ are ordered stochastically as their mean densities, that is: $\mu_1 \leq \mu_2$ if $\mu_1[\eta(0)] \leq \mu_2[\eta(0)]$. In particular, for every density $\rho \in [0, K]$, there is at most one $\nu \in (\mathcal{I} \cap \mathcal{S})_e$ with $\nu[\eta(0)] = \rho$; when it exists, we denote it by $\nu_\rho$.*

PROOF. The proof uses standard coupling arguments, which we briefly recall for self-containedness. We use couplings on the space $\widetilde{\mathbf{X}}$ and denote with a "tilde" quantities related to the coupled process. For any $\mu_1$ and $\mu_2$, elements of $(\mathcal{I} \cap \mathcal{S})_e$, there exists a coupling measure $\widetilde{\mu}$ with marginals $\mu_1$ and $\mu_2$ which belongs to the set $(\widetilde{\mathcal{I}} \cap \widetilde{\mathcal{S}})_e$ of invariant and translation invariant extremal measures for the coupled process (see [18] and [19]). By Theorem 1.9 in [8], $\widetilde{\mu}$ is supported on ordered configurations, that is, $(\eta, \xi)$ such that $\eta \leq \xi$ or $\xi \leq \eta$.

REMARK. To be precise we recall (as mentioned above in Section 2.1) that the difference between our model and the misanthrope's process studied in [8] is that for the latter additional assumptions are imposed on $b$ to ensure existence of product invariant measures. However these assumptions are *not* required for the proof of Theorem 1.9 in [8], which only uses assumptions (A1), (A3) and (A4), and extends in a natural way the coupling arguments introduced by Liggett [18] for the simple exclusion process. For details of the extension to misanthrope's process see Section 5 of [1] and [7].

It follows from extremality of $\widetilde{\mu}$ that the event $\{\eta \leq \xi\}$ has $\widetilde{\mu}$-probability 0 or 1 (cf. Chapter 8, Proposition 2.13 of [19]). This and $\mu_1[\eta(0)] \leq \mu_2[\xi(0)]$ imply $\widetilde{\mu}(\{\eta(0) \leq \xi(0)\}) = 1$. Indeed, $\widetilde{\mu}(\{\eta(0) > \xi(0)\}) > 0$ would imply $\widetilde{\mu}(\{\eta \leq \xi\}) = 0$, hence $\widetilde{\mu}(\{\xi \leq \eta\}) = 1$; but the latter combined with $\widetilde{\mu}(\{\eta(0) > \xi(0)\}) > 0$ in turn implies the contradiction $\mu_1[\eta(0)] = \widetilde{\mu}[\eta(0)] > \widetilde{\mu}[\xi(0)] = \mu_2[\xi(0)]$. Now, since $\widetilde{\mu}$ is translation invariant, $\widetilde{\mu}(\{\eta(0) \leq \xi(0)\}) = 1$ implies $\widetilde{\mu}(\{\eta(x) \leq \xi(x)\}) = 1$ for all $x \in \mathbb{Z}$. Hence $\widetilde{\mu}(\{\eta \leq \xi\}) = 1$, and thus $\mu_1 \leq \mu_2$. □

LEMMA 3.2. *$\mathcal{R}$ is a closed subset of $[0, K]$.*

PROOF. Since the empty and full configurations are unchanged by the dynamics, $\mathcal{R}$ contains at least 0 and $K$. Let $(\rho_n)_{n \geq 0}$ be an increasing sequence in $\mathcal{R}$, converging to $\rho$ (with obvious changes, what follows can be adapted to a decreasing subsequence). By Lemma 3.1, $(\nu_{\rho_n})_{n \geq 0}$ is an increasing sequence in $(\mathcal{I} \cap \mathcal{S})_e$, therefore converges (weakly) to some $\nu$; it is easy to see that $\nu$ belongs to $\mathcal{I} \cap \mathcal{S}$. Suppose $\rho \notin \mathcal{R}$, so $\nu$ is not extremal. Then $\nu = \int_{[0,K]} \nu_\alpha \gamma(d\alpha)$ for some measure $\gamma$ on $[0, K]$ such that $\gamma(\{\rho\}) = 0$, and

$$\nu[\eta(0)] = \lim_{n \to \infty} \nu_{\rho_n}[\eta(0)] = \rho = \int_{[0,K]} \alpha \gamma(d\alpha).$$

Then, there exists some $\rho' < \rho$ such that $\gamma([0, \rho']) > 0$, and by (iii) of Proposition 3.1,

(22) $$\nu\left\{\eta \in \mathbf{X} : \lim_{l \to \infty}(2l+1)^{-1} \sum_{x=-l}^{l} \eta(x) \leq \rho'\right\} > 0.$$

Let $n$ be such that $\rho_n > \rho'$, then because $\nu_{\rho_n} \leq \nu$, there exists a coupling measure $\widetilde{\nu}$, with marginals $\nu$ and $\nu_{\rho_n}$, such that $\widetilde{\nu}(\{(\eta, \xi) \in \widetilde{\mathbf{X}} : \xi \leq \eta\}) = 1$. This implies [using (iii) of Proposition 3.1]

$$\lim_{l \to \infty}(2l+1)^{-1} \sum_{x=-l}^{l} \eta(x) \geq \lim_{l \to \infty}(2l+1)^{-1} \sum_{x=-l}^{l} \xi(x) = \rho_n, \qquad \widetilde{\nu}\text{-a.s.}$$

Therefore $\nu(\{\eta \in \mathbf{X} : \lim_{l \to \infty}(2l+1)^{-1} \sum_{x=-l}^{l} \eta(x) \geq \rho_n\}) = 1$, which contradicts (22). Thus $\nu$ is extremal with density $\rho$. So $\rho$ belongs to $\mathcal{R}$, proving that $\mathcal{R}$ is closed. □

**4. The $\mathcal{R}$-valued Riemann problem.** In this section we derive local equilibrium for the so-called Riemann problem. This is a first step to obtain Theorem 2.2, for the particular initial condition (18) with $\lambda, \rho \in \mathcal{R}$. In other words we generalize [5], Section 2. We first define a corresponding initial measure $\nu_{\lambda,\rho}$ for the particle system which, as opposed to Andjel and Vares [2] and Bahadoran et al. [5], is generally not product. The measure $\nu_{\lambda,\rho}$ is defined in such a way as to enjoy the following properties:

(i) Negative (nonnegative) sites are distributed as under $\nu_\lambda$ ($\nu_\rho$);
(ii) $\tau_1 \nu_{\lambda,\rho} \geq \nu_{\lambda,\rho}$ ($\tau_1 \nu_{\lambda,\rho} \leq \nu_{\lambda,\rho}$) if $\lambda \leq \rho$ ($\lambda \geq \rho$);
(iii) $\nu_{\lambda,\rho}$ is stochastically increasing with respect to $\lambda$ and $\rho$.

In [5], $\nu_{\lambda,\rho}$ was product, and given explicitly by its marginals site per site. In the present setting we construct a (nonexplicit) $\nu_{\lambda,\rho}$ using the fact that if $\lambda \leq \rho$, then $\nu_\lambda \leq \nu_\rho$. Let $\widetilde{\mu}$ be a coupling measure on $\widetilde{\mathbf{X}}$ for $\nu_\lambda$ and $\nu_\rho$ such that $\widetilde{\mu}(\{(\eta,\xi): \eta \leq \xi\}) = 1$. We define the measure $\nu_{\lambda,\rho}$ as the distribution of the $\mathbf{X}$-valued random variable $\gamma = \gamma(\eta, \xi)$ defined on the coupling probability space $(\widetilde{\mathbf{X}}, \widetilde{\mu})$ by $\gamma(x) = \eta(x)$ if $x < 0$ and $\gamma(x) = \xi(x)$ if $x \geq 0$.

We shall prove in Section 4.2 below that the system starting from $\nu_{\lambda,\rho}$ is in local equilibrium with density profile given by the entropy solution to the Riemann problem with initial datum $u_0(\cdot)$. As a consequence we will obtain, in Corollary 4.2, the hydrodynamic limit for the $\mathcal{R}$-valued Riemann problem. As a first step, we prove in Section 4.1 that the $\mathcal{R}$-valued Riemann problem still admits the same variational characterization as in [5].

4.1. *A variational principle for Lipschitz-continuous flux.* The minimum regularity we must ensure for $G$ in order to extend the variational characterization of Bahadoran et al. [5] is Lipschitz continuity. We therefore prove the following:

LEMMA 4.1. *The flux $G$ defined in (17) is Lipschitz-continuous.*

PROOF. It is easy to see that there is a constant $C > 0$ and a finite set $S \subset \mathbb{Z}$ such that
$$|j(\eta) - j(\zeta)| \leq C \sum_{x \in S} |\eta(x) - \zeta(x)|.$$

For $\rho, r \in \mathcal{R}$, with $\rho \leq r$, we have $\nu_\rho \leq \nu_r$, so there is a coupling measure $\widetilde{\nu}_{\rho,r}$ for $\nu_\rho$ and $\nu_r$ under which $\eta \leq \zeta$ a.s.; thus
$$|G(\rho) - G(r)| \leq \widetilde{\nu}_{\rho,r}(|j(\eta) - j(\zeta)|) \leq C|S|(r - \rho).$$

Hence $G$ is uniformly Lipschitz inside $\mathcal{R}$, thus everywhere, since it is completed by linear interpolation and $\mathcal{R}$ is closed. $\square$

From now on we assume $\lambda < \rho$; adapting to the case $\lambda > \rho$ is straightforward, by replacing in the sequel lower with upper convex hulls, and minima/minimizers with maxima/maximizers. Consider $G_c$, the lower convex envelope of $G$ on $[\lambda, \rho]$. There exists a nondecreasing function $H_c$ (hence with left/right limits) such that $G_c$ has left/right hand derivative $H_c(\alpha \pm 0)$ at every $\alpha$; the notation $H_c(\alpha \pm 0)$ means left/right hand limit of $H_c$ at $\alpha$. The function $H_c$ is defined uniquely outside the at most countable set of nondifferentiability points of $G_c$

$$(23) \qquad \Theta = \{\alpha \in [\lambda, \rho] : H_c(\alpha - 0) < H_c(\alpha + 0)\}.$$

As will appear below, the particular choice of $H_c$ on $\Theta$ does not matter. Let $v_* = v_*(\lambda, \rho) = H_c(\lambda + 0)$ and $v^* = v^*(\lambda, \rho) = H_c(\rho - 0)$. Since $H_c$ is nondecreasing, there is a nondecreasing function $h_c$ on $[v_*, v^*]$ such that, for every $v \in [v_*, v^*]$,

$$(24) \qquad \begin{aligned} \alpha < h_c(v) &\implies H_c(\alpha) \leq v, \\ \alpha > h_c(v) &\implies H_c(\alpha) \geq v. \end{aligned}$$

Any such $h_c$ satisfies

$$(25) \qquad \begin{aligned} h_c(v - 0) &= \inf\{\alpha \in \mathbb{R} : H_c(\alpha) \geq v\} = \sup\{\alpha \in \mathbb{R} : H_c(\alpha) < v\}, \\ h_c(v + 0) &= \inf\{\alpha \in \mathbb{R} : H_c(\alpha) > v\} = \sup\{\alpha \in \mathbb{R} : H_c(\alpha) \leq v\}. \end{aligned}$$

It is not difficult to see that, anywhere in (25), $H_c(\alpha)$ may be replaced with $H_c(\alpha \pm 0)$. The following properties can be derived easily from (24) and (25):

1. Given $G$, $h_c$ is defined uniquely, and is continuous, outside the at most countable set

$$(26) \quad \begin{aligned} \Sigma_{\text{low}}(G) = \{v \in [v_*, v^*] : {}&G_c \text{ is differentiable with derivative } v \\ &\text{in a nonempty open subinterval of } [\lambda, \rho]\}. \end{aligned}$$

By "defined uniquely" we mean that for such $v$'s, there is a unique $h_c(v)$ satisfying (24), which does not depend on the choice of $H_c$ on $\Theta$.

2. Given $G$, $h_c(v \pm 0)$ is uniquely defined, that is, independent of the choice of $H_c$ on $\Theta$, for any $v \in [v_*, v^*]$. For $v \in \Sigma_{\text{low}}(G)$, $(h_c(v-0), h_c(v+0))$ is the maximal open interval over which $H_c$ has constant value $v$.

3. For every $\alpha \in \Theta$ and $v \in (H_c(\alpha-0), H_c(\alpha+0))$, $h_c(v)$ is uniquely defined and equal to $\alpha$.

In the sequel we extend $h_c$ outside $[v_*, v^*]$ in a natural way by setting

$$(27) \qquad h_c(v) = \lambda \quad \text{for } v < v_*, \qquad h_c(v) = \rho \quad \text{for } v > v^*.$$

PROPOSITION 4.1. *Let* $\lambda, \rho \in \mathcal{R}$, $\lambda < \rho$.

(i) *For every $v \in \mathbb{R} \setminus \Sigma_{\text{low}}(G)$, $G(\cdot) - v\cdot$ achieves its minimum over $[\lambda, \rho]$ at the unique point $h_c(v)$, and $u(x,t) := h_c(x/t)$ is the weak entropy solution to (7) with Riemann initial condition (18).*

(ii) *The previous minimum is unchanged if restricted to $[\lambda, \rho] \cap \mathcal{R}$. As a result, the Riemann entropy solution is a.e. $\mathcal{R}$-valued.*

This extends Proposition 2.1 of [5], where we assumed $G \in C^2(\mathbb{R})$. In the present setting part of the proof must be modified, because it explicitly used the regularity assumption. However the following lemma (Lemma 2.2 of [5]) carries over with the same proof.

LEMMA 4.2. *For $v \in \mathbb{R}$, $G(\cdot) - v\cdot$ and $G_c(\cdot) - v\cdot$ have the same global minimum value on $[\lambda, \rho]$; the set of global minimizers of $G_c(\cdot) - v\cdot$ on $[\lambda, \rho]$ is the interval $[h_c(v-0), h_c(v+0)]$; $h_c(v-0)$ and $h_c(v+0)$ are, respectively, the smallest and greatest global minimizer of $G(\cdot) - v\cdot$.*

As a consequence of Lemma 4.2 we obtain the following:

COROLLARY 4.1. (i) *If $v \in \Sigma_{\text{low}}(G)$, the graph of $G_c$ between $h_c(v-0)$ and $h_c(v+0)$ is a chord of the graph of $G$, with slope $v$, lying below the graph of $G$.*

(ii) *The graphs of $G$ and $G_c$ coincide at points $\alpha = h_c(v)$, where $v \in \mathbb{R} \setminus \Sigma_{\text{low}}(G)$.*

(iii) *Assume $\alpha \in [\lambda, \rho]$ and $G_c(\alpha) = G(\alpha)$. Then the restriction to $[\lambda, \alpha]$ (resp. $[\alpha, \rho]$) of $G_c$ is the lower convex envelope of the restriction of $G$ to $[\lambda, \alpha]$ (resp. $[\alpha, \rho]$).*

PROOF. (i) Let $v \in \Sigma_{\text{low}}(G)$. By Lemma 4.2, the points $(h_c(v \pm 0), G_c(h_c(v \pm 0)))$ lie on the graphs of $G$ and $G_c$. By convexity, $G_c(\cdot) - v\cdot$ is constant between the two global minima $h_c(v \pm 0)$, so the graph of $G_c$ between these points is linear with slope $v$. It is thus also a chord for $G$, and lies below $G$ because $G_c \leq G$.

(ii) Follows immediately from continuity of $h_c$ at $v$ and Lemma 4.2.

(iii) Denote by $G_c^1$ (resp. $G_c^2$) the lower convex hull of $G$ on $[\lambda, \alpha]$ (resp. $[\alpha, \rho]$). Then $G_c \leq G_c^1$ on $[\lambda, \alpha]$, since on this interval $G_c$ is convex and $G_c \leq G$. Similarly we have $G_c \leq G_c^2$ on $[\alpha, \rho]$. We now prove the converse inequalities. At endpoints the lower convex hull of a continuous function must coincide with the original graph (otherwise one could increase the convex hull in the neighborhood of an endpoint where it does not coincide, while keeping it convex and below $G$). Hence $G_c^1(\alpha) = G_c^2(\alpha) = G(\alpha) = G_c(\alpha)$. Denote by $H_c^i(\alpha \pm 0)$ the left/right hand derivative of $G_c^i$ at $\alpha$ for $i = 1, 2$. Because (as argued above) $G_c^1 \geq G_c$ on $[\lambda, \alpha]$, we must have $H_c^1(\alpha - 0) \leq H_c(\alpha - 0)$. Similarly $H_c^2(\alpha + 0) \geq H_c(\alpha + 0)$, and thus $H_c^1(\alpha - 0) \leq H_c^2(\alpha + 0)$. It follows that the function $G_{cc}$ obtained by concatenation of $G_c^1$ and $G_c^2$ is convex on $[\lambda, \rho]$. Since $G_{cc} \leq G$, we must have $G_{cc} \leq G_c$, which implies $G_c^1 \leq G_c$ (resp. $G_c^2 \leq G_c$) on $[\lambda, \alpha]$ (resp. $[\alpha, \rho]$). □

PROOF OF PROPOSITION 4.1. (i) We set $u(x,t) := h_c(x/t)$ for $x/t \notin \Sigma_{\text{low}}(G)$. Since $h_c$ is nondecreasing, it has bounded variation; it follows that $u$ has locally bounded space variation in the sense of (12). To prove that $u$ is the unique entropy solution given by Theorem 2.1 for Cauchy datum $u_0$, we check points (a), (b), (c) below: indeed the first two of them, combined with Proposition 2.2, imply that $u$ is an entropy solution:

(a) $u$ is a weak solution,
(b) jumps satisfy Oleĭnik's entropy condition (11),
(c) $u(\cdot, t)$ has initial datum $u_0$ in the sense of (9).

PROOF OF (a). Taking test functions, we see that $u(x,t) = h_c(x/t)$ is a weak solution of the hydrodynamic equation (7) if and only if $h_c(\cdot)$ satisfies

$$(28) \qquad \frac{d}{dv}[G(h_c(v)) - vh_c(v)] = -h_c(v)$$

weakly with respect to $v$. From Lemma 4.2,

$$G(h_c(v)) - vh_c(v) = \inf_{[\lambda,\rho]}[G(\cdot) - v \cdot] = \inf_{[\lambda,\rho]}[G_c(\cdot) - v \cdot] = -G_c^*(v)$$

for a.e. $v \in \mathbb{R}$, where $G_c^*$ is the Legendre–Fenchel transform of $G_c$ (see, e.g., [22]), therefore also a convex function, and $(G_c^*)' = h_c(v)$ almost surely with respect to Lebesgue measure. Using the absolute continuity of $G_c^*$ (which follows from convexity) (28) is obtained through integration by parts. □

PROOF OF (b). This follows from (i) of Corollary 4.1. □

PROOF OF (c). By (27), $u(\cdot, t)$ coincides with $u_0$ outside $[v_* t, v^* t]$. □

(ii) Assume that, for some $v \notin \Sigma_{\text{low}}(G)$, the unique minimum of $G(\cdot) - v \cdot$ is achieved at $\alpha \notin \mathcal{R}$. Then $\alpha \in (\lambda, \rho)$ and there exists an open interval $I \subset (\lambda, \rho)$ containing $\alpha$ such that $G$ is linear with slope $w \neq v$ on $I$ (recall that $G$ is linearly interpolated on $\mathcal{R}^c$). But then $G(\cdot) - v \cdot$ cannot achieve a minimum inside $I$. □

4.2. *Local equilibrium.* We now state the main result of this section:

PROPOSITION 4.2. *For every $v \in \mathbb{R} \setminus \Sigma_{\text{low}}(G)$, as $t \to \infty$,*

$$\nu_{\lambda,\rho}\tau_{[vt]}S(t) \Rightarrow \nu_{u(v,1)},$$

*where $u(x,t) = u(x/t,1)$ is the entropy solution to (7) constructed in Proposition 4.1 for the initial datum (18).*

A consequence is the following:

COROLLARY 4.2. *Assume $\eta_0^N \sim \nu_{\lambda,\rho}$ for every $N \in \mathbb{N}$. Then the sequence of processes $\eta_\cdot^N = (\eta_t^N, t \geq 0)$ has hydrodynamic limit $u(\cdot,\cdot)$ as above.*

REMARK. This is not an application of Landim [16], since we are no longer dealing with product measures.

PROOF OF COROLLARY 4.2. Pick a continuous test function $f$ with compact support in $\mathbb{R}$, and $l \in \mathbb{N}$. Then

$$\int f(x)\alpha^N(\eta, dx) = N^{-1}\sum_{x \in \mathbb{Z}} f(x/N)\eta(x)$$
$$= N^{-1}\sum_{x \in \mathbb{Z}} f(x/N)\eta^l(x) + o_N(1)$$
$$= \int f(x)\eta^l([Nx])\,dx + o_N(1),$$

where $\eta^l(x) = (2l+1)^{-1}\sum_{|y-x|\leq l}\eta(y)$, and $o_N(1)$ is a vanishing sequence depending on $f$ but not on $\eta$. Hence,

$$\begin{aligned}(29)\quad &\nu_{\lambda,\rho}\left(\left|\int f(x)\alpha^N(\eta_{Nt}, dx) - \int f(x)u(x/t,1)\,dx\right|\right)\\ &\leq \int f(x)\nu_{\lambda,\rho}\tau_{[Nx]}S(Nt)(|\eta^l(0) - u(x/t,1)|)\,dx + o_N(1).\end{aligned}$$

By Proposition 4.2 and (iii) of Proposition 3.1,

$$\begin{aligned}(30)\quad &\lim_{l\to\infty}\lim_{N\to\infty}\nu_{\lambda,\rho}\tau_{[Nx]}S(Nt)(|\eta^l(0) - u(x/t,1)|)\\ &= \lim_{l\to\infty}\nu_{u(x/t,1)}(|\eta^l(0) - u(x/t,1)|) = 0\end{aligned}$$

for a.e. $x \in \mathbb{R}$. The result follows from (29), (30) and Lebesgue's theorem. □

The proof of Proposition 4.2 is based on a scheme introduced by Andjel and Vares [2] for one-dimensional attractive processes with strictly concave (or convex) flux function $G$ and for which $(\mathcal{I} \cap \mathcal{S})_e$ consists of a continuous family of product measures indexed by their mean densities. In [5], we extended the proof up to nonconvex, nonconcave flux functions via a new variational formula. We generalize the latter proof, taking into account $\mathcal{R}$. Therefore we sketch the arguments that do not need modifications, and detail the others.

PROOF OF PROPOSITION 4.2. We assume $\lambda, \rho \in \mathcal{R}$ and $0 \leq \lambda < \rho$ (the case $\lambda > \rho$ is similar). The first step consists in showing that the limit in Proposition 4.2 holds in Cesáro sense. As in [5], we can transpose Lemmas 3.1 and 3.2 of [2] to our setting to get:

LEMMA 4.3.   (i) *Any sequence $T_n \to \infty$ has a subsequence $T_{n_m}$ for which there exists a dense set $D \subset \mathbb{R}$ such that for any $v \in D$*

$$\lim_{m \to \infty} \frac{1}{T_{n_m}} \int_0^{T_{n_m}} \nu_{\lambda,\rho} \tau_{[vt]} S(t) \, dt = \int \nu_\alpha \gamma_v(d\alpha) = \mu_v \in \mathcal{I} \cap \mathcal{S}$$

*for a probability measure $\gamma_v$ on $[\lambda, \rho] \cap \mathcal{R}$.*

(ii) *Moreover for any $u < v \in D$*

$$\lim_{m \to \infty} \nu_{\lambda,\rho} S(T_{n_m}) \left( \frac{1}{T_{n_m}} \sum_{[uT_{n_m}]}^{[vT_{n_m}]} \eta(x) \right) = F(u) - F(v)$$

*where $F(w) = \int [G(\alpha) - w\alpha] \gamma_w(d\alpha)$ for $w \in D$.*

The proof of (i) uses attractiveness, (i) and (ii) of Proposition 3.1, and the fact that the measure $\nu_{\lambda,\rho}$ satisfies $\nu_\lambda \leq \nu_{\lambda,\rho} \leq \tau_1 \nu_{\lambda,\rho} \leq \nu_\rho$; (ii) also uses the microscopic conservation equation (16). For these reasons the proofs of Andjel and Vares [2] carry over here without modification, as they made no use at that level of additional properties or assumptions (i.e., product invariant measures for all densities and strictly concave or convex flux $G$).

Next we can show, as in [5]:

LEMMA 4.4.   *For all $v \notin \Sigma_{\mathrm{low}}(G)$,*

(31)   $$\int_{[\lambda,\rho] \cap \mathcal{R}} [G(\alpha) - v\alpha] \gamma_v(d\alpha) \leq \inf_{\theta \in [\lambda,\rho] \cap \mathcal{R}} [G(\theta) - v\theta].$$

This is equation (9) from [5], except that here the minimum is restricted to $\mathcal{R}$. In fact we can exactly reproduce the proof given on pages 9–10 of [5]. Indeed, this proof uses only (a) attractiveness, (b) Lemma 4.3 above, (c) comparison of the process starting from $\nu_{\lambda,\rho}$ with the one starting from $\nu_{\theta,\rho}$ for intermediate densities $\theta \in [\lambda, \rho]$, and (d) comparison of the process starting from $\nu_{\theta,\rho}$ with the one starting from $\nu_\theta$. The only difference is that, in the present setting, we are constrained to choose $\theta \in \mathcal{R}$; this is why the restriction to $\mathcal{R}$ appears in (31). However from (ii) of Proposition 4.1, the infimum on the r.h.s. can be equivalently taken on the whole interval, and whenever it is uniquely achieved, this occurs in $\mathcal{R}$. Hence,

$$\gamma_v = \delta_{u(v,1)} \qquad \forall\, v \in \mathbb{R} \setminus \Sigma_{\text{low}}(G)$$

which concludes the first step. The second step is to replace the Cesáro limit by convergence in distribution. To this end, as in [5], we need Proposition 3.5 from Andjel and Vares [2]. This proposition uses only attractiveness, so it can be repeated here.

LEMMA 4.5. *Assume $\mu$ satisfies:*
(a) $\mu \leq \nu_\lambda$, (b) $\mu\tau_1 \geq \mu$, (c) *there exists $v_0$ finite so that*

$$\lim_{T\to\infty} \frac{1}{T} \int_0^T \mu\tau_{[vt]} S(t)\, dt = \nu_\lambda$$

*for all $v > v_0$. Then*

$$\lim_{t\to\infty} \mu\tau_{[vt]} S(t) = \nu_\lambda \qquad \text{for all } v > v_0.$$

This immediately implies

(32)
$$\lim_{t\to\infty} \nu_{\lambda,\rho}\tau_{[vt]} S(t) = \nu_\lambda \qquad \text{for all } v < v_*,$$
$$\lim_{t\to\infty} \nu_{\lambda,\rho}\tau_{[vt]} S(t) = \nu_\rho \qquad \text{for all } v > v^*.$$

Now consider $v \in [v_*, v^*] \setminus \Sigma_{\text{low}}(G)$. Set $\alpha = h_c(v)$. Let $(v_n)_n$ be a sequence such that $v_n < v$, $v_n \to v$ and $v_n \notin \Sigma_{\text{low}}(G)$ [the latter point is possible since $\Sigma_{\text{low}}(G)$ is countable]. Set $\alpha_n = h_c(v_n)$; then $\alpha_n \in \mathcal{R}$ by (ii) of Proposition 4.1. By monotonicity of $h_c$, $\alpha_n \leq \alpha$; by definition of $h_c$, $H_c(\alpha_n - 0) \leq v_n$ and thus $H_c(\alpha_n - 0) < v$; and by continuity of $h_c$ at $v$, $\alpha_n \to \alpha$.

Let $u^n(x,t) = u^n(x/t, 1)$ denote the entropy solution to (7) with initial datum $u_0^n(x) = \lambda \mathbf{1}_{\{x<0\}} + \alpha_n \mathbf{1}_{\{x\geq 0\}}$, and $G^n$ (resp. $G_c^n$) the restriction of $G$ (resp. $G_c$) to $[\lambda, \alpha_n]$. Since $v_n \notin \Sigma_{\text{low}}(G)$, by (ii), (iii) of Corollary 4.1, $G_c^n$ is the lower convex hull of $G^n$. Hence, for every $w > H_c(\alpha_n - 0) = v^*(\lambda, \alpha_n)$, we have, by (32), $\nu_{\lambda,\alpha_n}\tau_{[wt]} S(t) \to \nu_{\alpha_n}$ weakly for every such $w$. In particular we have $\nu_{\lambda,\alpha_n}\tau_{[vt]} S(t) \to \nu_{\alpha_n}$ weakly. By attractiveness, $\nu_{\lambda,\rho}\tau_{[vt]} S(t) \geq \nu_{\lambda,\alpha_n}\tau_{[vt]} S(t)$. Therefore any subsequential weak limit $\widetilde{\mu}_v$ of $\nu_{\lambda,\rho}\tau_{[vt]} S(t)$ satisfies $\widetilde{\mu}_v \geq \nu_{\alpha_n}$; since $\alpha_n \to \alpha$ we get $\widetilde{\mu}_v \geq \nu_\alpha$. To get the reverse inequality consider a sequence $v^n > v$, $v^n \to v$ with $v^n \notin \Sigma_{\text{low}}(G)$, set $\alpha^n = h_c(v^n)$ and consider $\nu_{\alpha^n,\rho}$. □

**5. General hydrodynamics.** In the previous section (Corollary 4.2) we established the hydrodynamic limit for certain initial measures associated with $\mathcal{R}$-valued Riemann initial profiles. In this section we prove that this implies Theorem 2.2, that is, the hydrodynamic limit for *any* initial sequence associated with *any* measurable initial density profile. To this end we adapt Section 3 of [5]. In that paper, we first proved Riemann hydrodynamics for *any* left- and right-hand densities. We then proceeded in the spirit of Glimm's scheme described in Section 2.4, by which one constructs general entropy solutions using only Riemann solutions. In our case we could only prove Riemann hydrodynamics for $\mathcal{R}$-valued left and right-hand densities, where $\mathcal{R}$ is a closed subset of $[0,K]$ outside which the macroscopic flux $G$ in (7) is linear. Therefore we refine the procedure in such a way as to reconstruct a general entropy solution by using only $\mathcal{R}$-valued Riemann solutions.

To extend hydrodynamics from Riemann to general initial profiles, we pointed out in [5] that the main property required for the particle system was what we called *macroscopic stability*. To state this property (Lemma 5.1 below) we need the following notation. Let $u(\cdot)$ and $v(\cdot)$ be two integrable $[0,K]$-valued density profiles on $\mathbb{R}$, and $\eta, \xi$ be two particle configurations in $\mathbf{X}$ with finitely many particles. Then we set

$$\Delta(u(\cdot), v(\cdot)) := \sup_{x \in \mathbb{R}} \left| \int_{-\infty}^{x} [u(y) - v(y)] \, dy \right|,$$

(33) $$\Delta^N(\eta, \xi) = N^{-1} \sup_{x \in \mathbb{Z}} \left| \sum_{y \leq x} [\eta(y) - \xi(y)] \right|,$$

$$\Delta^N(\eta, u(\cdot)) := N^{-1} \sup_{x \in \mathbb{Z}} \left| \sum_{y \leq x} \eta(y) - \int_{-\infty}^{x/N} u(y) \, dy \right|.$$

LEMMA 5.1. *Consider a sequence of coupled processes $(\eta^N_\cdot, \zeta^N_\cdot)$ (via the basic coupling) such that*

$$\sup_N N^{-1} \sum_{x \in \mathbb{Z}} [\eta^N_0(x) + \zeta^N_0(x)] < +\infty, \qquad a.s.$$

*Then, for every $t > 0$,*

(34) $$\Delta^N(\eta^N_{Nt}, \zeta^N_{Nt}) \leq \Delta^N(\eta^N_0, \zeta^N_0) + o_N(1)$$

*where $o_N(1)$ denotes a sequence of random variables converging to $0$ in probability.*

A consequence of (34) is that the hydrodynamic limit depends only on the density profile at time 0, and not on the underlying microscopic structure: that is, if a hydrodynamic limit holds for *some* initial sequence with given density profile, it still holds for *any* initial sequence with the same density profile. For the system defined by (2), (34) follows from the arguments of Bramson and Mountford [4], Proposition 3.1 (the latter is written for the asymmetric exclusion process; however the only properties involved are attractiveness, irreducibility and limited number of particles per site, so the proof can be adapted to our model). In the nearest-neighbor case [i.e., when $p(1) + p(-1) = 1$], (34) holds without even the $o_N(1)$ and is only a consequence of attractiveness, that is, assumption (A4). It follows also from a slightly more general result, that is, Proposition 2.1 of [3].

Another property we shall need is the *finite propagation* property for the particle system (see also Lemma 3.1 of [5]), which follows from the finite-range assumption on $p(\cdot)$. This is a microscopic counterpart to the finite propagation property for entropy solutions (Proposition 2.1).

LEMMA 5.2. *For any $x, y \in \mathbb{Z}$, any coupled process $(\eta_t, \zeta_t)_{t \geq 0}$ (via the basic coupling), and any $0 < t < (y - x)/(2V')$: if $\eta_0$ and $\zeta_0$ coincide on the site interval $[x, y]$, then $\eta_t$ and $\zeta_t$ coincide on the site interval $[x + V't, y - V't] \cap \mathbb{Z}$ with probability at least $1 - e^{-Ct}$, where $V'$ and $C$ are constants depending only on $b(\cdot, \cdot)$ and $p(\cdot)$.*

The proof of Theorem 2.2 from Corollary 4.2 can be decomposed into two steps. The first step is to prove the result for $\mathcal{R}$-valued entropy solutions. This step follows essentially the proof of Theorem 3.2 from [5], with some refinements due to $\mathcal{R}$. The second step, specific to this paper, is to show that an arbitrary entropy solution can be approximated by $\mathcal{R}$-valued entropy solutions. In the following, we shall recall the main steps followed in [5], but mainly insist on the differences that arise in the present setting.

5.1. *Hydrodynamics for $\mathcal{R}$-valued solutions.* In Theorem 3.2 of [5], we derived general hydrodynamics from Riemann hydrodynamics. We can repeat the same arguments here, when restricted to $\mathcal{R}$-valued entropy solutions. Before doing so, we need the following additional result:

LEMMA 5.3. *Assume the initial datum $u_0(\cdot)$ is a.e. $\mathcal{R}$-valued. Then, for every $t > 0$, the entropy solution $u(\cdot, t)$ to (7) is a.e. $\mathcal{R}$-valued.*

PROOF. The proof follows from Glimm's scheme. By Lemma 2.1 and (ii) of Proposition 4.1, the Glimm approximation (20)–(21) is a.e. $\mathcal{R}$-valued if $u_0$ is. Since $\mathcal{R}$ is closed, Theorem 2.3 implies the result.  □

We can now state the main result of this section:

PROPOSITION 5.1. *Theorem 2.2 holds when the initial density profile is a.e. $\mathcal{R}$-valued.*

The proof of this proposition, from the previously established $\mathcal{R}$-valued Riemann hydrodynamics (Corollary 4.2), is analogous to the proof of Theorem 3.2 on pages 16–17 of [5]. There we showed the following: for a system with finite initial configurations and compactly supported initial density profile, with the notations of Theorem 2.2 and (33), $\Delta^N(\eta^N_{Nt}, u(\cdot, t))$ cannot grow more than $o(\varepsilon)$ on a time interval of order $\varepsilon$. This implied the hydrodynamic limit for compactly supported initial profiles, and extension to general initial profiles followed from finite propagation arguments. The proof of this growth estimate for $\Delta^N$ relied on macroscopic stability, and Lemmas 3.5 and 3.6 of [5]. Lemma 3.6 allowed us to replace $u(\cdot, t)$, with error $o(\varepsilon)$, by a piecewise constant profile with step length at least $\varepsilon$. Lemma 3.5 showed that, starting from this approximate profile, the hydrodynamic limit at times of order $\varepsilon$ was indeed the entropy solution given by Lemma 2.1 as a succession of noninteracting Riemann waves. In the present case, Lemma 3.5 of Bahadoran et al. [5] is still true for $\mathcal{R}$-valued piecewise constant profiles:

LEMMA 5.4. *Assume $\varepsilon > 0$ and $u_0(\cdot)$ is an $\mathcal{R}$-valued step function with step lengths at least $\varepsilon$; then Theorem 2.2 holds up to time $\min(\varepsilon/(2V), \varepsilon/(2V'))$, where $V$ and $V'$ are the constants from Proposition 2.1 and Lemma 5.2.*

This lemma can be derived in the same way as Lemmas 3.4 and 3.5 of [5], from: (a) Riemann hydrodynamics (here Corollary 4.2 instead of Theorem 2.1 in the former paper), (b) macroscopic stability, (c) finite propagation property for particle systems and entropy solutions (Lemma 5.2 and Proposition 2.1 above). The idea is twofold: First, Riemann hydrodynamics were established in Corollary 4.2 for particular random initial configurations, but by macroscopic stability they can be extended to any initial random sequence with Riemann profile. Next, up to time $\min(\varepsilon/(2V), \varepsilon/(2V'))$, we have a succession of noninteracting Riemann problems, both on PDE and particle level: this is a consequence of Lemma 5.2 and Proposition 2.1. The only difference here compared to [5] is that, since point (a) is only established for $\mathcal{R}$-valued Riemann profiles, the conclusion also holds only for $\mathcal{R}$-valued step functions $u_0(\cdot)$.

At this stage a new difficulty appears compared to [5]: since the above lemma requires $\mathcal{R}$-valued piecewise constant profiles, we now need to approximate the entropy solution $u(\cdot, t)$—which we know is a.e. $\mathcal{R}$-valued by Lemma 5.3—with a piecewise constant profile also $\mathcal{R}$-valued. Hence Lemma 5.5 below is a refinement of Lemma 3.6 from [5], which requires a slightly different proof.

LEMMA 5.5. *Assume $u(\cdot)$ has compact support and finite variation, and is a.e. $\mathcal{R}$-valued. Let $\delta > 0$. Then, for $\varepsilon > 0$ small enough, there exists an approximation $u^{\varepsilon,\delta}$ of $u$ with the following properties: $u^{\varepsilon,\delta}$ is a piecewise*

EULER HYDRODYNAMICS OF ATTRACTIVE SYSTEMS 25*constant $\mathcal{R}$-valued function with compact support, step lengths at least $\varepsilon$, and $\Delta(u^{\varepsilon,\delta}, u) \leq \varepsilon\delta$.*

PROOF. Since $u$ has finite variation, it has left/right-hand limits $u(x \pm 0)$ at every $x \in \mathbb{R}$, and $u(x+0) = u(x-0)$ outside an at most countable set of points; note that $u(x \pm 0) \in \mathcal{R}$, because $\mathcal{R}$ is closed. Moreover, denoting by $\delta(x) := |u(x+0) - u(x-0)|$ the absolute jump at $x$, we have $\sum_{x \in \mathbb{R}} \delta(x) < +\infty$. Thus there is a finite set of points $J = \{x_1, \ldots, x_n\}$, with $x_1 < \cdots < x_n$, such that

$$\sum_{x \notin J} \delta(x) \leq \delta/4. \tag{35}$$

Let $x_0 < x_1$ and $x_{n+1} > x_n$ be such that the support of $u$ is contained in $[x_0, x_{n+1}]$. We can further divide each interval $I_k = (x_k, x_{k+1})$, for $k = 0, \ldots, n$, into open subintervals $I_{k,l}$ with lengths larger than $\varepsilon$ but smaller than $2\varepsilon$. For any interval $I \subset \mathbb{R}$, we set

$$\omega(I; \eta) = \operatorname{ess\,sup}\{|u(x) - u(y)| : x, y \in I, |x - y| \leq \eta\},$$
$$\omega(I) = \lim_{\eta \downarrow 0} \omega(I; \eta) = \sup_{x \in I} \delta(x).$$

In particular we have $\sup_{k=0,\ldots,n} \omega(I_k) \leq \delta/4$. Thus, for small enough $\varepsilon$,

$$\delta/4 + 2 \sup_{k=0,\ldots,n} \omega(I_k; 2\varepsilon) \leq \delta. \tag{36}$$

We shall define $u^{\varepsilon,\delta}$ as a piecewise constant function with value 0 outside $(x_0, x_{n+1})$ and constant value on each interval $I_{k,l}$. This value can be chosen as follows. Set $\rho_{k,l} = \operatorname{ess\,inf}_{I_{k,l}} u$, $\rho^{k,l} = \operatorname{ess\,sup}_{I_{k,l}} u$, $U(I_{k,l}) = \{u(x \pm 0), x \in I_{k,l}\} \subset \mathcal{R}$. Then

$$(37) \quad \forall \rho \in (\rho_{k,l}, \rho^{k,l}), \exists r \in (u(I_{k,l}) \cap \mathcal{R}) \cup U(I_{k,l}), \qquad |\rho - r| \leq \omega(I_{k,l})/2.$$

Hence there exists $r_{k,l} \in \mathcal{R}$ such that $|\overline{u}_{k,l} - r_{k,l}| \leq \omega(I_{k,l})/2$, where $\overline{u}_{k,l}$ denotes the mean value of $u$ on $I_{k,l}$. We set $u^{\varepsilon,\delta} \equiv r_{k,l}$ on $I_{k,l}$.

We claim that the function $u^{\varepsilon,\delta}$ thus defined satisfies the result of the lemma. Indeed, the integral $\int_{-\infty}^{x} [u(y) - u^{\varepsilon,\delta}(y)] \, dy$ has two contributions. The first one comes from integrating over intervals $I_{k,l}$ that do not contain $x$. By construction, this contribution is bounded in modulus by $2\varepsilon \sum_{k,l; x \notin I_{k,l}} \omega(I_{k,l})/2$. The second contribution is the integral over part of the unique $I_{k,l}$ containing $x$, that is, $\int_{x_{k,l}}^{x} [u(y) - u^{\varepsilon,\delta}(y)] \, dy$, where $x_{k,l}$ is the left-hand extremity of $I_{k,l}$. Since

$$|u(x) - r_{k,l}| \leq |u(x) - \overline{u}_{k,l}| + |\overline{u}_{k,l} - r_{k,l}|,$$

the latter contribution is bounded by $2\varepsilon\omega(I_{k,l}; 2\varepsilon) + 2\varepsilon\omega(I_{k,l})/2$. The result then follows from (35) and (36). $\square$

Proceeding from Lemmas 5.3, 5.4 and 5.5 above exactly as on pages 16–17 of [5] establishes Proposition 5.1.

5.2. *Hydrodynamics for arbitrary entropy solutions.* In order to approximate a general entropy solution with an $\mathcal{R}$-valued entropy solution, we use the following facts:

(i) an arbitrary initial datum with compact support can be $\Delta$-approximated with an $\mathcal{R}$-valued initial datum (even a $\{0,K\}$-valued initial datum is sufficient);

(ii) the entropy solution is $\Delta$-stable w.r.t. its initial datum (see [17]):

LEMMA 5.6. *Assume $u_0(\cdot)$ and $v_0(\cdot)$ have bounded support [then, by Proposition 2.1, so do the corresponding entropy solutions $u(\cdot,t)$ and $v(\cdot,t)$ at time t]. Then $\Delta(u(\cdot,t),v(\cdot,t)) \le \Delta(u_0(\cdot),v_0(\cdot))$.*

Note that this is a macroscopic analogue of Lemma 5.1. Equipped with the above results, we can now conclude the proof of Theorem 2.2 by deducing general hydrodynamics from $\mathcal{R}$-valued hydrodynamics. We shall use a technical result (whose proof is left to the reader):

LEMMA 5.7. *Assume $(\eta^N, N \in \mathbb{N}^*)$ is a sequence of $\mathbf{X}$-valued random variables such that:*

(i) *$\alpha^N(\eta^N, dx)$ [defined in (13)] converges in law along some subsequence to a random measure $\alpha(dx)$, necessarily of the form $\alpha(dx) = u(\cdot)\,dx$, where the random density profile $u(\cdot)$ is a.s. a.e. $[0,K]$-valued.*

(ii) *There exists some $a \in \mathbb{R}$ such that, with probability tending to 1 as $N \to \infty$, $\eta^N$ has no particle outside the site interval $[-aN, aN]$. Then, for every $v(\cdot): \mathbb{R} \to [0,K]$ with bounded support, $\Delta^N(\eta^N, v(\cdot))$ converges in law along the same subsequence to $\Delta(u(\cdot), v(\cdot))$.*

PROOF OF THEOREM 2.2. In the two steps below we use the following notation. $(\eta_0^N, N \in \mathbb{N}^*)$ denotes a sequence of random initial configurations with density profile $u_0(\cdot)$. The random measure $\alpha_t(dx) = u_t(\cdot)\,dx$ (where the random density profile $u_t(\cdot)$ is a.s. a.e. $[0,K]$-valued) denotes a subsequential limit in law of the sequence $\alpha^N(\eta_{Nt}^N, dx)$ as $N \to \infty$. Recall that $u(\cdot, t)$ denotes the entropy solution at time $t$ to (7) with Cauchy datum $u_0(\cdot)$. For $n \in \mathbb{N}$, $u_0^n(\cdot)$ denotes some approximation (to be defined) of $u_0(\cdot)$ as $n \to \infty$, $u^n(\cdot, \cdot)$ the corresponding entropy solution; $(\eta_0^{N,n}, N \in \mathbb{N}^*)$ a sequence (for fixed $n$) of random initial configurations with profile $u_0^n(\cdot)$ as $N \to \infty$, and $\eta_{Nt}^{N,n}$ the corresponding evolved configuration at time $Nt$. $\eta^N$ and $\eta^{N,n}$ are coupled via the basic coupling.

*Step* 1. We assume that, for some $a > 0$, $u_0(\cdot)$ vanishes outside $[-a, a]$ and $\eta_0^N$ has a.s. no particle outside $[-Na, Na]$. There is a sequence $(u_0^n(\cdot), n \in \mathbb{N})$ such that $u_0^n$ vanishes outside $[-a, a]$, is a.e. $\{0, K\}$-valued, and $\Delta(u_0^n, u_0) \to 0$ as $n \to \infty$. For instance,

$$u_0^n = K \sum_{i=0}^{n-1} \mathbf{1}_{[x_i, x_i + \alpha_i \delta)},$$

where $\delta = 2a/n$, $x_i = -a + i\delta$, and $\alpha_i = (K\delta)^{-1} \int_{x_i}^{x_{i+1}} u_0(x) \, dx$: indeed, since $u_0^n$ has the same integral as $u_0$ over each interval $[x_i, x_{i+1}]$, we have $\Delta(u_0^n, u_0) \leq K\delta = 2aK/n$. Given $u_0^n$, we can choose the sequence $(\eta_0^{N,n}, N \in \mathbb{N})$ with a.s. no particle outside $[-Na, Na]$.

Since $u_0^n(\cdot)$ is a.e. $\{0, K\}$-valued, $u^n(\cdot, \cdot)$ is a.e. $\mathcal{R}$-valued by Lemma 5.3. Section 5.1 established the hydrodynamic limit for $\mathcal{R}$-valued entropy solutions, hence $\alpha^N(\eta_{Nt}^{N,n}, dx)$ converges in probability to $u^n(\cdot, t) \, dx$ as $N \to \infty$ for each $n \in \mathbb{N}$. Furthermore, Lemma 5.2 implies that the sequences $\eta_{Nt}^{N,n}$ (for fixed $n$) and $\eta_{Nt}^N$ satisfy assumption (ii) of Lemma 5.7. Lemma 5.7 applied to both sequences, and the triangle inequality for $\Delta$, imply that along some subsequence $\Delta^N(\eta_{Nt}^{N,n}, \eta_{Nt}^N)$ converges in law to $\Delta(u^n(\cdot, t), u_t(\cdot))$ as $N \to \infty$. It follows from Lemma 5.1 that $\Delta(u^n(\cdot, t), u_t(\cdot)) \leq \Delta(u_0^n(\cdot), u_0(\cdot))$ a.s. with respect to the law of $u_t(\cdot)$. Hence $\Delta(u^n(\cdot, t), u_t(\cdot)) \to 0$ in probability as $n \to \infty$. By Lemma 5.6, $\Delta(u^n(\cdot, t), u(\cdot, t)) \to 0$ as $n \to \infty$; hence $u_t(\cdot) = u(\cdot, t)$ a.s.

*Step* 2. Now we make no further assumption on $u_0(\cdot)$ and $\eta_0^N$. For $n \in \mathbb{N}$ we set $u_0^n = u_0 \mathbf{1}_{[-n,n]}$ and $\eta_0^{N,n}(x) = \eta_0^N(x) \mathbf{1}_{[-Nn, Nn]}(x)$. By Lemma 5.2, $\eta_{Nt}^N$ and $\eta_{Nt}^{N,n}$ coincide on the site interval $[N(-n + V't), N(n - V't)] \cap \mathbb{Z}$ with probability tending to 1 as $N \to \infty$. This and step one applied to $\eta_{Nt}^{N,n}$ implies that $u_t(\cdot)$ coincides with $u^n(\cdot, t)$ a.s. on $[-n + V't, n - V't]$. Proposition 2.1 implies $u^n(\cdot, t) = u(\cdot, t)$ a.e. on $[-n + Vt, n - Vt]$. Letting $n \to \infty$ yields $u_t(\cdot) = u(\cdot, t)$ a.s. $\square$

## APPENDIX A: PROOF OF PROPOSITION 2.2

Vol'pert [26] considers the class BV of real-valued functions $u(x, t)$ on $\mathbb{R} \times \mathbb{R}^{+*}$ defined as follows: $u \in BV$ if and only if, for every compact subset $K$ of $\mathbb{R} \times \mathbb{R}^{+*}$, there exists a constant $C(K)$ such that

$$\left| \int_K u(x,t) \, \partial_t \varphi(x, t) \, dx \, dt \right| \leq C(K) \sup_K |\varphi|,$$

$$\left| \int_K u(x,t) \, \partial_x \varphi(x, t) \, dx \, dt \right| \leq C(K) \sup_K |\varphi|,$$

for every smooth test function $\varphi$ with support contained in $K$. Equivalently, $u \in BV$ if and only if its gradient in the sense of distributions is a locally finite vector measure. To summarize the results of Vol'pert [26] on the structure of BV functions, we need the following notions. According to Definition 2.1 of [6], we say that $u$ has an *approximate jump discontinuity* at $y = (x, t)$ in the direction of vector $\nu \neq 0$, if there exist $u^- \neq u^+$ such that

$$(38) \qquad \lim_{r \to 0} r^{-2} \int_{z=(x,t): |z| \leq r} |u(y+z) - U_{u^-, u^+, \nu}(z)| \, dx \, dt = 0,$$

where

$$U_{u^-, u^+, \nu}(z) = \begin{cases} u^+, & \text{if } \nu \cdot z > 0, \\ u^-, & \text{if } \nu \cdot z < 0, \end{cases}$$

and $\nu \cdot z$ denotes scalar product in $\mathbb{R}^2$. For $\nu = (1, -v)$, the meaning is that $u$ has a local discontinuity at $y$ traveling with speed $v$, with value $u^-$ on the left and $u^+$ on the right. On the other hand, if (38) holds with $u^- = u^+$, then it actually holds for any $\nu \neq 0$, and we say $y$ is an *approximate continuity point* for $u$. The main results are: (i) Outside a set of one-dimensional Hausdorff measure 0, every point is either a point of approximate continuity or a point of approximate jump discontinuity. (ii) The set of approximate jump points is of Hausdorff dimension 1 and has locally finite one-dimensional Hausdorff measure.

If $u$ has locally bounded space variation, its spatial derivative $\partial_x u$ in distribution sense is a locally finite measure. Following (7) and Lipschitz continuity of $G$, its time derivative $\partial_t u$ is also a locally finite measure; hence $u \in BV$. One technical problem arises from the coexistence of two a priori different notions of limits for $u$. Indeed, (i) on the one hand we may view $u(\cdot, t)$ as a function with locally bounded space variation for all $t$, hence it has one-sided limits at every $x \in \mathbb{R}$ for every fixed $t$; (ii) on the other hand, $u$ is a BV function of the two variables, and as such it has approximate jump discontinuities as defined in (38). Limits in the sense (i) do not have a normal space–time vector (or equivalently a local velocity) attached to them, so they may not be limits in the sense (ii). Conversely, limits in the sense (ii) are not defined for all $t$ and $x$, but only in space–time average, so they may not be limits in the sense (i). However (a spatially localized version of) Theorem 2.6 of [6] shows that these two notions coincide, provided we can prove the additional estimate

$$(39) \qquad \int_x^y |u(z, t) - u(z, s)| \, dz \leq M|t - s|$$

for any $0 < s < t$ and $x, y \in \mathbb{R}$, where the constant $M$ is uniform when $(s, t, x, y)$ varies in a bounded set. Equation (39) is a consequence of (7) and Lipschitz continuity of $G$. Theorem 2.6 of Bressan [6] states that: (i) up to a set of one-dimensional Hausdorff measure 0, the jump set of $u$ as a BV function of $(x, t)$ in the sense of Vol'pert [26] coincides with the set $\{(x, t) \in \mathbb{R} \times \mathbb{R}^{+*}; u(x - 0, t) \neq (x + 0, t)\}$, (ii) the local normal is almost nowhere (w.r.t. one-dimensional Hausdorff measure) parallel to the time axis, (iii) for a jump point $(x,t)$, let $u^+$ (resp. $u^-$) denote the approximate limit on the positive (resp. negative) side of the normal chosen with positive $x$-coordinate. Then, a.e. w.r.t. one-dimensional Hausdorff measure, we have $u^\pm = u(x \pm 0, t)$.

Equation (8) is satisfied by all entropy-flux pairs, if it merely holds for the so-called Kružkov entropy-flux pairs $(\phi_c, \psi_c)$ defined for $c \in [0, K]$ by

$$\phi_c(u) := |u - c|, \qquad \psi_c(u) = \text{sgn}(u - c)[G(u) - G(c)]$$

(see, e.g., Chapter 2 of [25]). The first theorem on page 256 of [26] states that a weak solution $u \in BV$ of (7) satisfies (8) for $(\phi, \psi) = (\phi_c, \psi_c)$ if and only if

$$(40) \qquad [\phi_c(u^+) - \phi_c(u^-)]\nu_t + [\psi_c(u^+) - \psi_c(u^-)]\nu_x \leq 0$$

holds for a.e. jump point of $u$ w.r.t. one-dimensional Hausdorff measure, where $u^+$ (resp. $u^-$) is the approximate limit of $u$ on the positive (resp. negative) side of the normal vector $\nu = (\nu_x, \nu_t)$. [Note: this is an elementary result if $u$ is piecewise continuously differentiable, because then the l.h.s. of (40) is the trace of the l.h.s. of (8) in distribution sense along the shock line, when $\phi = \phi_c$ and $\psi = \psi_c$.] Equation (40) does not depend on the choice of $\nu$, but by (ii) above we may choose $\nu$ with positive $x$-coordinate; then, by (iii), $u^\pm = u(x \pm 0, t)$. Since the l.h.s. of (40) is continuous w.r.t. $c$, (40) is satisfied simultaneously for all $c$ outside a single exceptional set. We conclude by noticing that the following are equivalent for a pair $(u^-, u^+)$:

(a) (40) holds for every $c \in [0, K]$,
(b) $(u^-, u^+)$ satisfies Oleinik's entropy condition (11).

To see this, first take $c = 0$ to obtain that $-\nu_t/\nu_x$ is the slope of the chord of $G$ between $u^-$ and $u^+$; then let $c$ vary between $u^-$ and $u^+$.

## APPENDIX B: PROOF OF THEOREM 2.3

In the sequel we incorporate dependence on the discretization step $\Delta x$ and on the sampling sequence $(a) = (a_k)_{k \in \mathbb{Z}^+} \in (-1, 1)^{\mathbb{Z}^+}$ in the Glimm approximation of $u$, denoting it by $\tilde{u}^{\Delta x, (a)}$. We will denote by $\nu$ the product uniform measure on $(-1, 1)^{\mathbb{Z}^+}$.

For $m \in \mathbb{N}$, let $u_{m,0}$ denote an approximation of $u_0$ such that: (i) $u_{m,0}$ has locally bounded variation, (ii) $u_{m,0}$ converges to $u_0$ in $L^1_{\text{loc}}(\mathbb{R})$ as $m \to \infty$. We denote by $u_m$ the entropy solution to (7) with initial datum $u_{m,0}$, and by $\tilde{u}_m^{\Delta x, (a)}$ the Glimm approximation of $u_m$. We know that Theorem 2.3 is already true when the initial datum has locally bounded variation, hence for $u_{m,0}$. Thus, by diagonal extraction, there is a sequence $\delta_n \downarrow 0$ (as $n \to \infty$) and a subset $\mathcal{A}$ of $(-1, 1)^{\mathbb{Z}^+}$ with $\nu(\mathcal{A}) = 1$, such that

$$(41) \quad \forall (a) \in \mathcal{A}, m \in \mathbb{N} \qquad \lim_{\Delta x = \delta_n \downarrow 0} \tilde{u}_m^{\Delta x, (a)} = u_m \qquad \text{in } L^1_{\text{loc}}(\mathbb{R} \times \mathbb{R}^{+*}).$$

On the other hand, by (10), $u_m$ converges to $u$ in $L^1_{\text{loc}}(\mathbb{R} \times \mathbb{R}^{+*})$ as $m \to \infty$. The idea is now to compare the Glimm approximations $\tilde{u}^{\Delta x,(a)}$ and $\tilde{u}^{\Delta x,(a)}_m$. Set $V'' := V + R^{-1}$, with $V$ the propagation speed in (10) and $R$ the ratio defined in (19). We claim that, for every $x < y$ and $t < (y-x)/(2V'')$,

$$(42) \quad \mathbb{E}_\nu \left[ \int_{x+V''t}^{y-V''t} |\tilde{u}^{\Delta x,(a)}_m(z,t) - \tilde{u}^{\Delta x,(a)}(z,t)| \, dz \right] \leq \int_x^y |u_0(z) - u_{m,0}(z)| \, dz,$$

where $\mathbb{E}_\nu$ means that $(a)$ is integrated w.r.t. $\nu$. The meaning of (42) is that the Glimm approximations "nearly" satisfy the same contraction–propagation property (10) (with a different propagation speed $V''$) as the real entropy solutions. However "nearly" means that this is not true for a fixed choice of the sampling sequence $(a)$, but in average with respect to this sequence. Before proving (42), let us conclude the proof of Theorem 2.3: for any bounded time interval $I \subset (0, (y-x)/(2V''))$,

$$(43) \quad \begin{aligned} \mathbb{E}_\nu & \left[ \int_I \int_{x+V''t}^{y-V''t} |u(z,t) - \tilde{u}^{\Delta x,(a)}(z,t)| \, dz \, dt \right] \\ & \leq \int_I \int_{x+V''t}^{y-V''t} |u(z,t) - u_m(z,t)| \, dz \, dt \\ & \quad + \mathbb{E}_\nu \left[ \int_I \int_{x+V''t}^{y-V''t} |u_m(z,t) - \tilde{u}^{\Delta x,(a)}_m(z,t)| \, dz \, dt \right] \\ & \quad + \mathbb{E}_\nu \left[ \int_I \int_{x+V''t}^{y-V''t} |\tilde{u}^{\Delta x,(a)}_m(z,t) - \tilde{u}^{\Delta x,(a)}(z,t)| \, dz \, dt \right]. \end{aligned}$$

The first term on the r.h.s. of (43) vanishes as $m \to \infty$ by (10); the second one vanishes as $\Delta x = \delta_n \downarrow 0$ for any fixed $m$, by (41) and dominated convergence; and the third one vanishes as $\Delta x \to 0$ and then $m \to \infty$ by (42). We have thus proved that, for any $[x,y] \subset \mathbb{R}$ and bounded time interval $I \subset (0, (y-x)/(2V''))$, the limit

$$(44) \quad \lim_{\Delta x \to 0} \int_I \int_{x+V''t}^{y-V''t} |u(z,t) - \tilde{u}^{\Delta x,(a)}(z,t)| \, dz \, dt = 0$$

holds in $L^1(\nu)$ along the subsequence $\Delta x = \delta_n$. Hence, taking growing countable families of intervals $[x,y]$ and $I$, and using diagonal extraction, we can construct a further subset $\mathcal{A}'$ of $(-1,1)^{\mathbb{Z}^+}$ with $\nu(\mathcal{A}') = 1$, and a further subsequence $(h_n)$ of $(\delta_n)$, such that (44) holds along $\Delta x = h_n$ for every $(a) \in \mathcal{A}'$, $[x,y] \in \mathbb{R}$ and $I \subset (0, (y-x)/(2V''))$. This concludes the proof of Theorem 2.3.

We now proceed to prove (42). Let us introduce the following notations: $(S_t)_{t \geq 0}$ denotes the evolution semigroup defined by the entropy solution to (7), that is, $S_t u_0$ is the entropy solution at time $t$ when the initial datum is $u_0$; for $a_k \in (-1, 1)$ and $\Delta x > 0$, $T^{a_k, \Delta x}$ denotes the approximation operator defined by $\tilde{u}_k^+ = T^{a_k, \Delta x} \tilde{u}_k^-$ in (21). Thus, for $t \in [t_k, t_{k+1})$,

$$(45) \qquad \tilde{u}^{\Delta x, (a)}(\cdot, t) = S_{t-t_k} T^{a_k, \Delta x} S_{\Delta t} T^{a_{k-1}, \Delta x} \cdots S_{\Delta t} T^{a_0, \Delta x} u_0.$$

It follows from definition of $T^{a_k, \Delta x}$ that

$$(46) \quad \mathbb{E}_\nu \left[ \int_x^y |T^{a_k, \Delta x} v(z) - T^{a_k, \Delta x} w(z)| \, dz \right] \leq \int_{x-\Delta x}^{y+\Delta x} |v(z) - w(z)| \, dz$$

for every measurable $[0, K]$-valued functions $v(\cdot)$ and $w(\cdot)$ on $\mathbb{R}$, and every subinterval $[x, y] \subset \mathbb{R}$. On the other hand, by (10), we also have

$$(47) \qquad \int_x^y |S_{\Delta t} v(z) - S_{\Delta t} w(z)| \, dz \leq \int_{x-V\Delta t}^{y+V\Delta t} |v(z) - w(z)| \, dz;$$

(42) follows from (45), (46) and (47).

**Acknowledgments.** We thank the referees for their careful reading of the paper and their suggestions. We thank Tom Mountford for interesting discussions. We thank IHP for hospitality during the semesters *Limites hydrodynamiques* and *Geometry and statistics of random growth*. H. Guiol thanks Université de Rouen and EPFL, and E. Saada thanks E.P.F.L. for hospitality.

C. BAHADORAN
LABORATOIRE DE MATHÉMATIQUES
UNIV. CLERMONT-FERRAND 2
63177 AUBIÈRE
FRANCE
E-MAIL: bahadora@math.univ-bpclermont.fr

H. GUIOL
TIMC–TIMB
INP GRENOBLE
FACULTÉ DE MÉDECINE
38706 LA TRONCHE CEDEX
FRANCE
E-MAIL: Herve.Guiol@imag.fr

K. RAVISHANKAR
SUNY NEW PALTZ
COLLEGE AT NEW PALTZ
1 HAWK DRIVE
NEW PALTZ, NEW YORK 12561
USA
E-MAIL: ravishak@newpaltz.edu

E. SAADA
CNRS, UMR 6085
UNIVERITÉ DE ROUEN
SITE DU MADRILLET
AVENUE DE L'UNIVERSITÉ
BP. 12, 76801
SAINT ETIENNE DU ROUVRAY
FRANCE
E-MAIL: Ellen.Saada@univ_rouen.fr